\documentclass[oneside,a4paper]{amsart}
\usepackage{poly-min-surf}

\begin{document}

\title{Area minimizing polyhedral surfaces are saddle}
\author{Anton Petrunin}
\address{Anton Petrunin, 
Math. Dept., 
PSU,
University Park, 
PA 16802, USA.}
\address{petrunin@math.psu.edu}
\maketitle


\section*{Preface}

At age fifteen
I had to solve the following problem:

\begin{thm}{Problem}
Consider all quadrangles $\square axby$ in the plane with fixed sides $|a-x|$, $|a-y|$, $|b-x|$ and $|b-y|$.
Note that the value 
$$\alpha=\measuredangle axb+\measuredangle ayb$$ 
describes the quadrangle $\square axby$ up to congruence;
let $A(\alpha)$ be the area of the quadrangle for given $\alpha$.

Show that $A(\alpha)$ increases in $\alpha$ for $\alpha\le \pi$ 
and decreases in $\alpha$ for $\alpha\ge\pi$.%
\footnote{In particular the area of a quadrilateral with fixed side lengths is maximal when it is inscribed into a circle.}
\end{thm}

The problem was not especially hard, beautiful or interesting.
But a voice in my head said ``one day it will be useful'' ---
a strange warning that turned out to be true.

Ten years later 
I was finishing graduate school.
I was trying to prove something about minimal surfaces in Hadamard spaces 
(not important what it is).
As I simplified the problem further and further, 
I eventually saw the problem above.

It proved what I wanted to and made me happy for few days.
Later, I generalized the statement yet further
and it ended up in my paper \cite{petrunin}.
The technique used in the original proof turned out to be redundant.
On the other hand, 
the argument was simple and beautiful, 
so I decided that it was worth sharing.

In addition I notice recent closely related activity, 
see for example \cite{bobenko},
but as far as I can see the idea below has not been noticed.

\section*{Introduction}

The following question is a simplified version of the one mentioned in the preface;
still it contains all of its interesting features.

Let $\DD$ be a simplicial complex homeomorphic to the disc in the plane
(think of a convex polygon with fixed triangulation).
A piecewise linear map $F$ from $\DD$ to the Euclidean space
 will be called a \emph{polyhedral disc};
that is, a map
$F\:\DD\to \RR^3$ is a polyhedral disc if the restriction of $F$ to any triangle of $\DD$ is linear.
Intuitively the polyhedral disc is a disc in $\RR^3$ glued from triangles with possible self-intersections.

With slight abuse of notation,
we make no distinction between 
vertices, edges, and triangles of  $\DD$
and their $F$-images.

We say that a vertex or an edge of $F$ is \emph{interior} if does not lie in the boundary $\partial\DD$.

The \emph{area} of a polyhedral disc $F$ is defined as the sum of the 
areas of all its triangles.

A polyhedral disc is called \emph{saddle} if one can not cut a hat from it by a plane.
More precisely,
we say that a plane $\Pi$ cuts an edge $[ab]$ 
if the endpoints $a$ and $b$ lie on opposite sides of $\Pi$.
Then the polyhedral disc $F$ is called saddle if for any interior vertex $a$
of $F$
there is no plane 
which cuts each edge coming from $a$.

Fix a positive integer $n$.
Consider a class of polyhedral discs $F$
with the same boundary curve $F(\partial\DD)$
and with the total number of triangles 
at most $n$.
A disc $F$ is called \emph{area minimizing} 
if it has the minimal area in this class.

\begin{thm}{Theorem}
Any area minimizing polyhedral disc in Euclidean space is saddle.
\end{thm}

Before we get into the proof, 
let us discuss an example. 
Assume we change the definition of area minimizing polyhedral disc 
a bit;
instead of giving the upper bound for the number of triangles,
we fix one triangulation.
In this case, 
the conclusion of theorem does not hold.

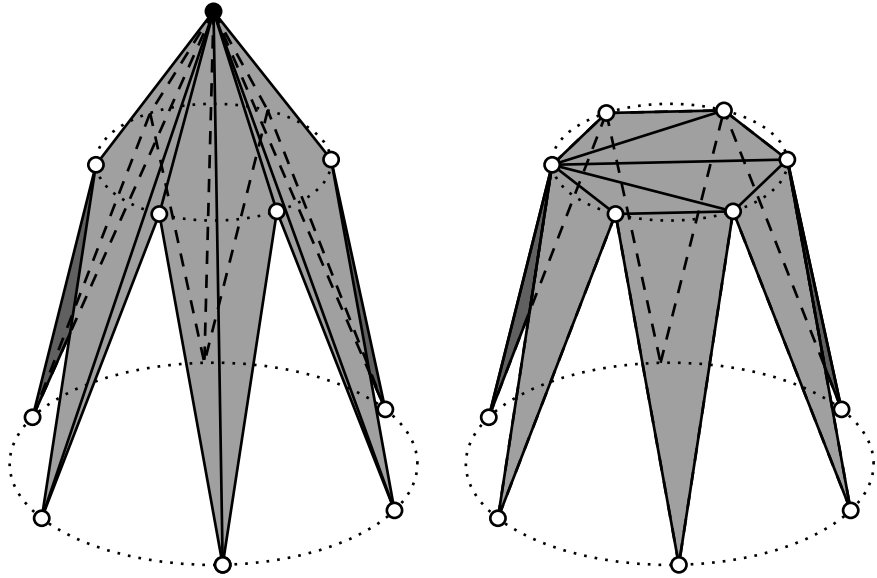
\begin{center}
\begin{lpic}[t(-0mm),b(-0mm),r(0mm),l(0mm)]{pics/tent(1)}
\end{lpic}
\end{center}

A counterexample is shown on the left.
This tent forms a polyhedral disc made from 12 triangles;
it is mirror symmetric in each vertical plane which contains one white 
and the black vertex on the top.
The black vertex is the only interior vertex of the disc;
it can be cut by a plane from the rest, so this disc is not saddle.

The disc on the left minimizes the area for the given triangulation.
The disc on the right has smaller area 
and a smaller number of triangles (there are 10 of them).
In fact, the disc on the right is area minimizing for $n=10$;
it has to be saddle according to the theorem,
but it is also saddle by a trivial reason ---
it has no interior vertices.

\section*{Proof}

Let $F$ be an area minimizing polyhedral disc.

Without loss of generality,
we can assume that one can not reduce the number of triangles in $F$ 
while keeping the area the same.
In this case,
$F$ satisfies the so-called \emph{no triangle condition};
i.e.,
if three vertices of $F$, say $a$, $b$ and $c$,  
are pairwise joined by edges 
then $\triangle abc$ is a triangle of $F$.

Indeed, if this is not the case,
exchange the domain bounded by these three edges by $\triangle abc$;
this procedure will not increase the area of $F$ and it will drop
the number of triangles in the triangulation.

\parbf{Claim.}
\textit{We can assume that the sum of all 4 angles which are adjacent to any interior edge of $F$ is at least $\pi$.
Moreover if this sum is $\pi$, 
then the two adjacent triangles together form a flat convex quadrilateral.}

\medskip

Assume the contrary;
i.e., there is an edge $[ab]$ in $F$ 
with two adjacent triangles $\triangle abx$ and $\triangle aby$
such that 
\[\measuredangle abx+\measuredangle aby+\measuredangle bax+\measuredangle bay<\pi.\leqno\bigstar\]
Let us cut these two triangles from $F$
and glue $\triangle axy$ and $\triangle bxy$ instead.
This way we get a new polyhedral disc, say $H$,
with a new triangulation.

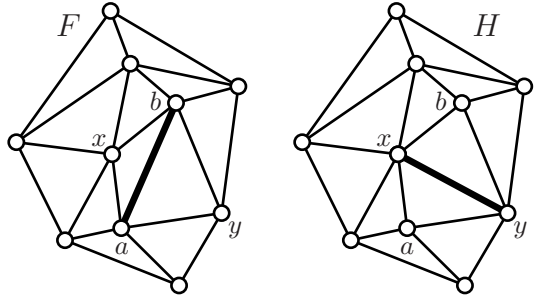
\begin{center}
\begin{lpic}[t(-0mm),b(-0mm),r(0mm),l(0mm)]{pics/flip(1)}
\lbl{8,37;{\Large $F$}}
\lbl{63,37;{\Large $H$}}
\lbl[t]{15,8.2;$a$}
\lbl[r]{20.2,27;$b$}
\lbl[rb]{13,21;$x$}
\lbl[tl]{29,11;$y$}
\lbl[t]{52.5,8.2;$a$}
\lbl[r]{57.7,27;$b$}
\lbl[rb]{50.5,21;$x$}
\lbl[tl]{66.5,11;$y$}
\end{lpic}
\end{center}

The construction of $H$ from $F$
will be called \emph{flip of the edge $[ab]$}.
Note that performing the flip we will always get a genuine triangulation;
this follows since the original triangulation 
satisfies the no triangle condition.

Let us show that 
\[\area F\ge \area H.\leqno\spadesuit\]
To do this, we construct two quadrilaterals 
$\square a'x'b'y'$ and $\square a''x''b''y''$ in the plane
such that the diagonal $[a'b']$ divides $\square a'x'b'y'$,
the diagonal $[x''y'']$ divides $\square a''x''b''y''$
and
\begin{align*}
\triangle abx&\cong\triangle a'b'x',
&
\triangle axy&\cong\triangle a''x''y'',
\\
\triangle aby&\cong\triangle a'b'y',
&
\triangle bxy&\cong\triangle b''x''y''.
\end{align*}

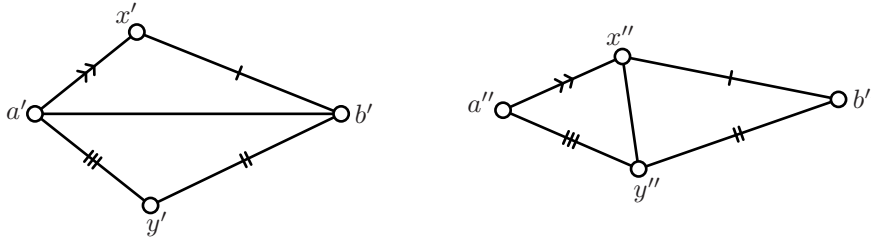
\begin{center}
\begin{lpic}[t(1mm),b(1mm),r(0mm),l(0mm)]{pics/two-trig(1)}
\lbl[r]{2,15.5;$a'$}
\lbl[l]{45,15;$b'$}
\lbl[b]{15,27;$x'$}
\lbl[t]{19,2;$y'$}
\lbl[r]{63.5,16.5;$a''$}
\lbl[l]{110.5,17;$b''$}
\lbl[b]{80,24;$x''$}
\lbl[t]{83.5,6.5;$y''$}
\end{lpic}
\end{center}

Note that
\[
\begin{aligned}
\measuredangle x'a'y'+\measuredangle x'b'y'
&=\measuredangle xab+\measuredangle yab+\measuredangle xba+\measuredangle yba\ge
\\
&\ge \measuredangle xay+\measuredangle xby =
\\
&=\measuredangle x''a''y''+\measuredangle x''b''y''.
\end{aligned}
\leqno\clubsuit\]
Applying $\bigstar$ and the Problem,
we get
\[\area (\square a'x'b'y')\ge \area (\square a''x''b''y''),\]
or equivalently,
\[\area(\triangle abx)+\area(\triangle aby)
\ge \area(\triangle axy)+\area(\triangle bxy).\leqno\diamondsuit\]
Hence $\spadesuit$ follows.

Note that we have equality in $\diamondsuit$
if and only if we have equality in $\clubsuit$.
Further, equality in $\clubsuit$ holds if and only if the quadrilateral  
$\square axby$ is flat and convex.
Therefore, 
if the disc $F$ is in general position; 
i.e., no 4 vertices of the disc lie in one plane, 
then the Claim follows.

Further, we can assume that the triangulation of $F$ is chosen in such a way that there is an approximation of $F$ by discs in general position such that no flip decreases its area.
Hence the Claim follows in the general case.

Now assume the disc $F$ is not saddle.
In this case we can move one of its interior vertices, say $a$, 
so that all the edges coming from $a$ become shorter.
To do this, choose a plane which cuts each edge coming from $a$,
and move $a$ toward the plane along the segment perpendicular to the plane, say with unit speed.
Let us denote by $a(t)$ the position of $a$ after time $t$
and let $F_t$ be the obtained polyhedral disc.

In general this deformation may not decrease the area.
However it does decrease the area for the discs which satisfy the statement in the Claim.

Indeed, note that the area of $F$ is completely determined by the triangulation and the lengths of its interior edges.
Assume $\ell_1(t),\dots,\ell_k(t)$ are the lengths of the  edges coming from $a(t)$.
Then 
\[\area F_t=A(\ell_1(t),\dots,\ell_k(t)).\]
Applying the Problem again, 
we get that
\[\frac{\partial A}{\partial\ell_i}\ge 0 \leqno\heartsuit\]
for each $i$.
Thus, $t\mapsto\area F_t$ is decreasing for small $t$ 
if for at least one $i$ the inequality $\heartsuit$ is strict.

Finally note that if  for each $i$
we get equality in 
$\heartsuit$, 
then the sum of 4 adjacent angles at each edge from $a$ is exactly $\pi$.
Therefore, from the second statement in the Claim, 
all the edges coming from $a$ lie in one plane.
These edges can not point into  a fixed open half-plane, simply because an angle of a triangle can not be bigger than $\pi$.
In particular, there is no plane which cuts the edges coming from $a$,
a contradiction.\qeds

\parbf{Acknowledgment.}
I wish to thank  my adviser, Stephanie Alexander.

\end{document}